# ASYMPTOTIC PROPERTIES OF A FAMILY OF SOLUTIONS OF THE PAINLEVÉ EQUATION $P_{VI}$

OVIDIU COSTIN AND RODICA D. COSTIN

October 24, 2018

## 1. Introduction and setting

In analyzing the question whether nonlinear equations can define new functions with good global properties, Fuchs had the idea that a crucial feature now known as the *Painlevé property* (PP) is the absence of movable (meaning their position is solution-dependent) essential singularities, primarily branch-points, see [8]. First order equations were classified with respect to the PP by Fuchs, Briot and Bouquet, and Painlevé by 1888, and it was concluded that they give rise to no new functions. Painlevé and Gambier took this analysis to second order, looking for all equations of the form $u'' = F(u', u, z)$, with $F$ rational in $u'$, algebraic in $u$, and analytic in $z$, having the PP [17, 18]. They found some fifty types with this property and succeeded to solve all but six of them in terms of previously known functions. The remaining six types are now known as the Painlevé equations. Beginning in the 1980's, almost a century after their discovery, these equations were related to linear problems (and thereby solved) by various methods including the powerful techniques of isomonodromic deformation and reduction to Riemann-Hilbert problems [3], [4], [7], [15]. The solutions of the six Painlevé equations play a fundamental role in many areas of pure and applied mathematics due to their integrability properties. In particular, there are numerous physical applications of the Painlevé $P_{VI}$ equation (for references see e.g. [6]) among which we mention the problem of construction of self-dual Bianchi-type IX Einstein metrics, [2, 5, 16, 22] the classification of the solutions of the Witten-Dijkgraaf-Verlinde-Verlinde equation (WDVV) in 2D-topological field theories, and probability theory, especially random matrix theory (see e.g. [23], [24]). The connection between determinants and Painlevé equations was established in the early 70's (see [19], [20] and references therein). The two point correlation functions for holonomic fields on the Poincaré disk are shown to be expressible in terms of in terms of $P_{VI}$ [21].

A three parameter family of solutions of the Painlevé equation $P_{VI}$ arises in the context of random matrix theory in a recent work of Borodin and Deift [1].

The asymptotic behavior of the solutions of the Painlevé equations is of utmost importance. The main purpose of this paper is to characterize a family of solutions of $P_{VI}$ for large argument, relevant to the study [1]

In the $\sigma$-form these solutions satisfy (see [9]—eq. (C.61), with $\nu_1 = \nu_2$):

$$u'[u''t(t-1)]^2 + \left[2u'(tu'-u) - u'^2 - \nu_1^2\nu_3\nu_4\right]^2 = (u'+\nu_1^2)^2(u'+\nu_3^2)(u'+\nu_4^2) \tag{1.1}$$

where $\Re\nu_1 > 0$.

Eq. (1.1) admits the exact solution

$$u(t) = -{\nu_1}^2 t + \frac{1}{2}\left({\nu_1}^2 + \nu_3\,\nu_4\right)$$

and a one parameter family of solutions with the behavior for large $t$

$$u(t) = -{\nu_1}^2 t + \frac{1}{2}\left({\nu_1}^2 + \nu_3\,\nu_4\right) + Ct^{-2\nu_1} + O(t^{-2\nu_1-1}) \tag{1.2}$$

**Proposition 1.** *For any $C \in \mathbb{C}$, eq. (1.1) has a unique solution satisfying*

$$u(t) = -{\nu_1}^2 t + \frac{1}{2}\left({\nu_1}^2 + \nu_3\,\nu_4\right) + Ct^{-2\nu_1} + o\left(t^{-2\nu_1}\right)$$

*for $t \to \infty$ in any fixed sector $\mathcal{S}$.*

## 2. Proof of Proposition 1

2.1. **Notation.** Denote by $\mathcal{F}$ the set of functions of the form $a(t) = f(t^{-1}, t^{-\nu_1})$ where $f$ is analytic at $(0,0)$. Note that functions in $\mathcal{F}$ are bounded for $t$ large enough.

Let $C \in \mathbb{C}$ and denote

$$u(t) = -{\nu_1}^2 t + \frac{1}{2}\left({\nu_1}^2 + \nu_3\,\nu_4\right) + Ct^{-2\nu_1} + \Delta(t) \tag{2.3}$$

We only consider $\Delta$ such that

$$\Delta(t) = o\left(t^{-2\nu_1}\right) \tag{2.4}$$

for large $t$ in $\mathcal{S}$. Then we also have, as is easy to see by the Cauchy formula,

$$\Delta'(t) = o\left(t^{-2\nu_1-1}\right) \quad , \quad \Delta''(t) = o\left(t^{-2\nu_1-2}\right) \tag{2.5}$$

for large $t$ in the sector.

2.2. **Equation for the remainder $\Delta$.** Substituting (2.3) in (1.1) we get the equation

$$T_2\,\left(\Delta''\right)^2 + T_1\,\Delta'' + T_0 = 0 \tag{2.6}$$

where $T_j$ depend on $t, \Delta$ and $\Delta'$ and have the form

$$\begin{aligned} T_2 &= a_2(t) - \Delta' \\ T_1 &= t^{-2-2\nu_1} b_1(t) + t^{-2-2\nu_1} a_1(t)\Delta' \\ T_0 &= t^{-5-4\nu_1} c_2(t) + t^{-4-2\nu_1} c_0(t)\Delta + t^{-3-2\nu_1} c_1(t)\Delta' + t^{-2}P \end{aligned} \tag{2.7}$$

where $a_j, b_j, c_j \in \mathcal{F}$, and

$$P = d_1(t)t^{-2}\Delta'\Delta^2 + t^{-1}d_2(t)\Delta'\Delta + t^{-2}d_3(t)\Delta^2 \tag{2.8}$$

$$+(\Delta')^2\left[q_1(t) + t^{-1}q_2(t)\Delta + t^{-2}q_3(t)\Delta^2 + q_4(t)\Delta' + t^{-1}q_5(t)\Delta\Delta' + q_6(t)(\Delta')^2\right]$$

with $d_j, q_j \in \mathcal{F}$ (the Appendix contains exact formulae). We also have

$$a_2(t) = {\nu_1}^2 + O\left(t^{-1}\right)$$

(2.9)
$$c_1(t) = 8{\nu_1}^4 C(2\nu_1 + 1) + O\left(t^{-1}\right)$$

$$c_0(t) = -8{\nu_1}^4 C(2\nu_1 + 1) + O\left(t^{-1}\right)$$

We write eq. (2.6) in normal form, solved for $\Delta''$, and we separate the dominant terms.

**Lemma 2.** *The function* $\Delta$ *satisfies the equation*

(2.10)
$$\Delta'' + 2\,\nu_1 t^{-1}\Delta' - 2\,\nu_1 t^{-2}\Delta = R$$

*where* $R$ *depends on* $t$, $\Delta$ *and* $\Delta'$, *and gathers smaller terms*

(2.11)
$$R = t^{-3-2\nu_1}\tilde{c}_4 + F_{1s} + \tilde{R}_2 - \tau$$

*where the terms are given by (3.30), (3.32), (3.33), (3.34) and (2.15).*

*Proof of Lemma 2.*
From (2.6) we have

(2.12)
$$\Delta'' = \frac{-T_1 \pm \sqrt{T_1^2 - 4T_0T_2}}{2T_2}$$

The minus choice in (2.12) is not consistent with (2.5). Indeed, since $T_0T_2/T_2^2 = o(1)$ we have

$$\frac{-T_1 - \sqrt{T_1^2 - 4T_0T_2}}{2T_2} = -\frac{1}{2}\,\frac{T_1}{T_2}\,\left(2 + o(1)\right)$$

which is of order $t^{-2\nu_1-2}$, hence is not $o\left(t^{-2\nu_1-2}\right)$.
Thus

(2.13)
$$\Delta'' = \frac{-T_1 + \sqrt{T_1^2 - 4T_0T_2}}{2T_2} = \frac{-2T_0}{T_1}\,\frac{1}{1 + \sqrt{1 - \frac{4T_0T_2}{T_1^2}}} \equiv F(t, \Delta, \Delta')$$

To separate the dominant linear part of equation (2.13) we rewrite $F$ as

(2.14)
$$F(t, \Delta, \Delta') = -\frac{T_0}{T_1} - \tau$$

where

(2.15)
$$\tau = 4\left(\frac{T_0}{T_1}\right)^2 \frac{T_2}{T_1}\frac{1}{\left(1 + \sqrt{1 - \frac{4T_0T_2}{T_1^2}}\right)^2}$$

A direct calculation of $T_0/T_1$ yields (2.10) (see §3.5 for details). ∎

2.3. **Integral equations for $\Delta$ and $\Delta'$.** The left hand-side of (2.10) has the solutions $t$ and $t^{-2\nu_1}$, hence equation (2.10) can be written in integral form as

$$\Delta(t) \,=\, \frac{1}{2\nu_1+1}\,\left[t\int_\infty^t R(s)\,ds \,-\, t^{-2\nu_1}\,\int_\infty^t s^{1+2\nu_1}\,R(s)\,ds\right] \tag{2.16}$$

Denote

$$\Delta_1 = \Delta \ , \ \Delta_2 = \Delta' \tag{2.17}$$

Equation (2.16) becomes the system of first order integral equations for $(\Delta_1, \Delta_2)$

(2.18)

$$\Delta_1(t) \,=\, \tfrac{1}{2\nu_1+1}\,\left[t\int_\infty^t R(s)\,ds \,-\, t^{-2\nu_1}\,\int_\infty^t s^{1+2\nu_1}\,R(s)\,ds\right] \,\equiv\, J_1(\Delta_1,\Delta_2)$$

$$\Delta_2(t) \,=\, \tfrac{1}{2\nu_1+1}\,\left[\int_\infty^t R(s)\,ds \,+ 2\nu_1\, t^{-1-2\nu_1}\,\int_\infty^t s^{1+2\nu_1}\,R(s)\,ds\right] \,\equiv\, J_2(\Delta_1,\Delta_2)$$

2.4. **Existence and uniqueness of $\Delta = O\left(t^{-1-2\nu_1}\right)$.** Consider the domain

$$D = \{t \in \mathbb{C}\,;\, |t| > \rho\ ,\ \arg t \in (A,B)\,\}$$

where $\rho$ will be chosen large enough and $A < B < A + 2\pi$. (Sectors of larger angles can be considered on the Riemann surface above $\mathbb{C}\setminus 0$).

Let $\mathcal{B}$ be the Banach space of pairs $\Delta = (\Delta_1, \Delta_2)$ of analytic function on $D$, continuous on $\overline{D}$ with

$$\|\Delta\| \equiv \max\left\{\sup_{t\in D} t^{1+2\nu_1}|\Delta_1(t)|\ ,\ \sup_{t\in D} t^{2+2\nu_1}|\Delta_2(t)|\right\} \,<\, \infty$$

We will show that the integral operator $\mathbf{J} = (J_1, J_2)$ defined by (2.18) applies a ball of $\mathcal{B}$ into itself and is a contraction there. This implies that (2.18) has a unique solution in $\mathcal{B}$.

2.4.1. $\mathbf{J}$ *applies a ball of $\mathcal{B}$ into itself.* Let $\mathcal{B}_M$ be the ball of elements of $\Delta \in \mathcal{B}$ of norm at most $M$: we have $|\Delta_1| \le M\,t^{-1-2\nu_1}$ and $|\Delta_2| \le M\,t^{-2-2\nu_1}$.

From (3.32)

$$|F_{1s}| \,\le\, const\,|t|^{-4}\,\left(|t|^{-2\nu_1}|\Delta| \,+\, |t|^{-2\nu_1-1}|\Delta'|\right) \tag{2.19}$$

hence for $(\Delta, \Delta') \in \mathcal{B}_M$

$$|F_{1s}| \,\le\, const\,|t|^{-4-2\nu_1}\,M \tag{2.20}$$

To estimate $\tilde{R}_2$ from (3.33) we note that for $(\Delta, \Delta') \in \mathcal{B}_M$ we have $|S_{1,2}| \le |t|^{-1}$ (see (3.29) for notations) and we have

(2.21)

$$|\tilde{R}_2| \le const\ |t|^{-4-2\nu_1}\,M^2\,\left[1 + M|t|^{-2-2\nu_1} + \left(M|t|^{-2-2\nu_1}\right)^2\right]\left(1 - M|t|^{-2-2\nu_1}\right)^{-1}$$

From (2.11), (2.20), (2.21), (3.38), (3.40), (3.41) we get

$$|R| \le K\ |t|^{-3-2\nu_1}\ \left(1 + t^{-1}\Phi(M,t)\right) \tag{2.22}$$

where

$$\Phi(M,t) = M + M^2\left[1 \,+\, M|t|^{-2-2\nu_1} \,+\, \left(M|t|^{-2-2\nu_1}\right)^2\right]\left(1 \,-\, M|t|^{-2-2\nu_1}\right)^{-1}$$
$$+ \left[1 + M + M^2t^{-1} + (M^3 + M^4)t^{-3-2\nu_1}\right]^2 (1 + Mt^{-2-2\nu_1})\,(1 - Mt^{-2-2\nu_1})^{-3}$$

and $K$ is independent of $M$.

For $|t| > \rho$ (2.22) implies

$$|R| \leq K\ |t|^{-3-2\nu_1}\ \left(1 + \rho^{-1}\Phi(M,\rho)\right) \tag{2.23}$$

therefore, from (2.18) we get

$$|J_1(\Delta_1,\Delta_2)| \leq K'\, t^{-1-2\nu_1}\ \left(1 + \rho^{-1}\Phi(M,\rho)\right)$$

and

$$|J_2(\Delta_1,\Delta_2)| \leq K'\, t^{-2-2\nu_1}\ \left(1 + \rho^{-1}\Phi(M,\rho)\right)$$

Choosing $M > K'$ and then $\rho_0$ large enough it follows that for any $\rho > \rho_0$ the operator $\mathbf{J}$ applies the ball $\mathcal{B}_M$ into itself.

2.4.2. $\mathbf{J}$ *is a contraction on* $\mathcal{B}_M$. The parameter $M$ is now fixed by §2.4.1 (hence the constants in the estimates of the present section may depend on $M$); $\rho$ will be chosen large enough.

Let $\Delta^{[1]}, \Delta^{[2]}$ be two elements in $\mathcal{B}_M$.

From (3.32) we see that

$$|F_{1s}(\Delta^{[1]}) - F_{1s}(\Delta^{[2]})| \leq const\ t^{-4-2\nu_1}\ \|\Delta^{[1]} - \Delta^{[2]}\|$$

and, from (3.33) we get

$$|\tilde{R}_2(\Delta^{[1]}) - \tilde{R}_2(\Delta^{[2]})| \leq const\ t^{-4-2\nu_1}\ \|\Delta^{[1]} - \Delta^{[2]}\| \tag{2.24}$$

(see §3.6 for details).

We also have

$$|\tau(\Delta^{[1]}) - \tau(\Delta^{[2]})| \leq const\ t^{-4-2\nu_1}\ \|\Delta^{[1]} - \Delta^{[2]}\| \tag{2.25}$$

The details are in §3.7.

Then from (2.18)

$$\|\mathbf{J}(\Delta^{[1]}) - \mathbf{J}(\Delta^{[2]})\| \leq K\ t^{-1}\ \|\Delta^{[1]} - \Delta^{[2]}\|$$

which shows that $\mathbf{J}$ is a contraction on $\mathcal{B}_M$ if $\rho$ is large enough.

Then (2.18) has a unique solution in $\mathcal{B}_M$.

2.5. **Uniqueness of** $\Delta = o\left(t^{-2\nu_1}\right)$**.** Let $\Delta$ be a solution of (2.10) satisfying $\Delta = o(t^{-2\nu_1})$ for large $t$ in a sector. We now show that, in fact, $\Delta = O(t^{-1-2\nu_1})$ which completes the proof of Proposition 1.

Note that we have $\Delta' = o(t^{-1-2\nu_1})$.

For any $\epsilon > 0$ there exists $\rho > 0$ such that

$$|\Delta| \leq \epsilon\, |t|^{-2\nu_1} \quad , \quad |\Delta'| \leq \epsilon\, |t|^{-1-2\nu_1}$$

for $|t| > \rho$ in the sector.

From (3.31) we get

$$|F_{1s}| \leq const\ \epsilon\, |t|^{-3-2\nu_1} \tag{2.26}$$

From (3.33) and the estimates of §3.8.1 we get

$$|\tilde{R}_2| \leq const\ \epsilon\ \left(|t^{-2}\Delta| + |t^{-1}\Delta'|\right)\ +\ const\ \epsilon^2\, |t^{-3}T^{-4}| \tag{2.27}$$

Using (3.38), (3.42) and the estimates detailed in §3.8.2 we get that also $\tau$ has an upper bound of the form of the RHS of (2.27).

Then from (2.11), (2.26), (2.27) it follows that

$$|R| \le const \ \left( |t^{-3-2\nu_1}| + \epsilon\, |t^{-2}\Delta| + \epsilon\, |t^{-1}\Delta'| + \epsilon^2\, |t^{-3}T^{-4}| \right) \tag{2.28}$$

This shows that (in the notations (2.17)) the integrals in (2.18) are convergent, so the solution $\Delta$ of (1.1) satisfies (2.18).

Denote

$$\|(\Delta_1, \Delta_2)\| = \sup_t \{|\Delta_1|, |t\Delta_2|\}$$

Using (2.28) in (2.18) we get

$$|\Delta_1| \le const \ \left( |t^{-1-2\nu_1}| + \epsilon \|(\Delta_1, \Delta_2)\| \right)$$

and

$$|\Delta_2| \le const \ \left( |t^{-2-2\nu_1}| + \epsilon\, |t|^{-1} \, \|(\Delta_1, \Delta_2)\| \right)$$

hence

$$\|(\Delta_1, \Delta_2)\| \le const \ \sup |t^{-1-2\nu_1}| \quad \text{for } |t| > \rho$$

so that $\Delta = O(t^{-1-2\nu_1})$ and the proof of Proposition 1 is complete.

## 3. Appendix

### 3.1. Notation used.

$$T = t^{\nu_1}$$

### 3.2. The expression of $T_2$.

$$T_2 = -\Delta'(t) - \frac{\nu_1 \left( 2\,\nu_1\, t - t^2 \nu_1 - \nu_1 - 2\, t^{-1-2\,\nu_1} C + 4\, C t^{-2\,\nu_1} - 2\, t^{1-2\,\nu_1} C \right)}{(t-1)^2}$$

### 3.3. The expression of $T_1$.

$$(t-1)^2 T_1 = 4C\nu_1 \left(2\nu_1 + 1\right) \left( -t^{-2\nu_1} - t^{-2-2\nu_1} + 2t^{-1-2\nu_1} \right) \Delta'(t) + 4C{\nu_1}^2 \left(2\nu_1 + 1\right)$$
$$\left[ t^{-2-2\nu_1}\nu_1 + t^{-2\nu_1}\nu_1 - 2t^{-1-2\nu_1}\nu_1 + 2t^{-1-4\nu_1} C - 4t^{-2-4\nu_1} C + 2t^{-3-4\nu_1} C \right]$$

### 3.4. The expression of $T_0$.

$$T_0 = L_{00} + L_{01}\Delta + L_{02}\Delta' + R_{em}$$

where

$$\begin{aligned} L_{00} = & -4\nu_1^3 C^2 \Big[ \left( 4{\nu_1}^3 + 6{\nu_1}^2 + 4\nu_1\nu_3\nu_4 + 2\nu_3\nu_4 + 2\nu_1 \right) t \\ & \quad - 4{\nu_1}^2 - \nu_1 + \nu_1{\nu_4}^2 - 4{\nu_1}^3 - 2\nu_1\nu_3\nu_4 + \nu_1{\nu_3}^2 \Big] \frac{1}{t^4 (t-1)^2 T^4} \\ & - 8C^3{\nu_1}^2 \left( \left( 4{\nu_1}^2 + 4{\nu_1}^3 + \nu_1 \right) t^2 + \left( 2\nu_1 + 4{\nu_1}^2 + 4\nu_1\nu_3\nu_4 + 2\nu_3\nu_4 \right) t \right. \\ & \quad \left. -4{\nu_1}^2 - \nu_1 + \nu_1{\nu_4}^2 - 4{\nu_1}^3 - 2\nu_1\nu_3\nu_4 + \nu_1{\nu_3}^2 \right) \frac{1}{t^5 (t-1)^2 T^6} \\ & - 16 \frac{(1+2\nu_1)\, C^4 {\nu_1}^2 \left( (1+2\nu_1)\, t - 2\nu_1 \right)}{t^5 (t-1)^2 T^8} \end{aligned}$$

and

$$L_{01} = -8\frac{C{\nu_1}^3\left(\left(\nu_1+2{\nu_1}^2\right)t+\nu_3\nu_4-{\nu_1}^2\right)}{t^3\left(t-1\right)^2T^2}$$
$$-16\frac{{\nu_1}^2C^2\left(\left(2\nu_1+4{\nu_1}^2\right)t-2{\nu_1}^2+\nu_3\nu_4\right)}{t^4\left(t-1\right)^2T^4}-32\frac{C^3{\nu_1}^2\left(\left(1+2\nu_1\right)t-\nu_1\right)}{t^5\left(t-1\right)^2T^6}$$

Also

$$L_{02} = \frac{4}{t^3\left(t-1\right)^2T^2}\nu_1^2C\Big[\left(2{\nu_1}^2+4\nu_1^3\right)t^2+\left(-{\nu_1}^2+4\nu_1\nu_3\nu_4+\nu_3\nu_4-4{\nu_1}^3\right)t$$
$$+\nu_1{\nu_3}^2-2\nu_1\nu_3\nu_4+\nu_1{\nu_4}^2\Big]$$
$$+\frac{4}{t^4\left(t-1\right)^2T^4}\nu_1C^2\left(\left(12{\nu_1}^2+\nu_1+20{\nu_1}^3\right)t^2+\left(12\nu_1\nu_3\nu_4-16\nu_1^3+4\nu_3\nu_4+2\nu_1\right)t\right.$$
$$\left.+3\nu_1\nu_3^2-4{\nu_1}^3-4{\nu_1}^2+3\nu_1{\nu_4}^2-\nu_1-6\nu_1\nu_3\nu_4\right)$$
$$+\frac{16}{t^4\left(t-1\right)^2T^6}\nu_1C^3\left(\left(6\nu_1+1+8{\nu_1}^2\right)t-8{\nu_1}^2-3\nu_1\right)$$

Finally

$$R_{em} = -4\frac{{\Delta'}^4}{\left(t-1\right)t}+\frac{\left(8t-4\right){\Delta'}^3\Delta}{t^2\left(t-1\right)^2}$$
$$+\left(\frac{8t^2{\nu_1}^2+4\nu_3\nu_4t-2\nu_3\nu_4+{\nu_4}^2-8{\nu_1}^2t+{\nu_3}^2}{t^2\left(t-1\right)^2}+4\frac{C\left(-8\nu_1+2t+8\nu_1t-1\right)}{t^2\left(t-1\right)^2T^2}\right){\Delta'}^3$$
$$-4\frac{{\Delta'}^2\Delta^2}{t^2\left(t-1\right)^2}+\left(-2\frac{C\left(24\nu_1t-12\nu_1+4t\right)}{t^3\left(t-1\right)^2T^2}-\frac{16{\nu_1}^2t+4\nu_3\nu_4-8{\nu_1}^2}{t^2\left(t-1\right)^2}\right){\Delta'}^2\Delta$$
$$+\left(-2\frac{1}{t^3\left(t-1\right)^2T^2}C\left(2\nu_3\nu_4t+8t^2{\nu_1}^2-6\nu_1\nu_3\nu_4+12\nu_1\nu_3\nu_4t-24\nu_1^3t+3\nu_1{\nu_3}^2\right.\right.$$
$$\left.+24{\nu_1}^3t^2+3\nu_1{\nu_4}^2\ -4{\nu_1}^2t\right)-4\frac{C^2\left(t+12\nu_1t-24{\nu_1}^2-6\nu_1+24{\nu_1}^2t\right)}{t^3\left(t-1\right)^2T^4}$$
$$\left.-\frac{-4{\nu_1}^4t-2{\nu_1}^2\nu_3\nu_4+4{\nu_1}^2\nu_3\nu_4t+{\nu_1}^2{\nu_4}^2+{\nu_3}^2{\nu_1}^2+4{\nu_1}^4t^2}{t^2\left(t-1\right)^2}\right){\Delta'}^2$$
$$+\left(8\frac{{\nu_1}^2}{t^2\left(t-1\right)^2}+16\frac{C\nu_1}{t^3\left(t-1\right)^2T^2}\right)\Delta'\Delta^2$$
$$+\left(16\frac{C^2\nu_1\left(2t-3\nu_1+6\nu_1t\right)}{\left(t-1\right)^2t^4T^4}+4\frac{{\nu_1}^2\left(-{\nu_1}^2+2{\nu_1}^2t+\nu_3\nu_4\right)}{t^2\left(t-1\right)^2}\right.$$
$$\left.+16\frac{C\nu_1\left(4{\nu_1}^2t+\nu_3\nu_4-2{\nu_1}^2+\nu_1t\right)}{t^3\left(t-1\right)^2T^2}\right)\Delta'\Delta-4\frac{{\nu_1}^4\Delta^2}{t^2\left(t-1\right)^2}-16\frac{\Delta^2{\nu_1}^3C}{t^3\left(t-1\right)^2T^2}$$
$$-16\left(\frac{\Delta\nu_1C}{(t-1)t^2T^2}\right)^2$$

3.5. **Splitting of terms of $T_0/T_1$.** We introduce the notations:

$$S_1 = t^{2\nu_1}\Delta \quad , \quad S_2 = t^{1+2\nu_1}\Delta' \tag{3.29}$$

Note that in the assumptions of Proposition 1 we have $S_1, S_2 = o(1)$, and for $(\Delta,\Delta') \in \mathcal{B}$ we have $S_1, S_2 = O(t^{-1})$.

Separating the terms of $T_0/T_1$ by degree and dominance we have

$$\frac{T_0}{T_1} = F_0 + F_{1d} + F_{1s} + \tilde{R}_2$$

where

$$F_0 = \frac{c_2}{b_1 t^3 T^2} \equiv \tilde{c}_4 t^{-3-2\nu_1} \tag{3.30}$$

The linear terms are

$$F_{1d} + F_{1s} = \frac{c_0\Delta}{b_1 t^2} + \left(\frac{c_1}{b_1} + \frac{c_2 a_1}{t^2 T^2 {b_1}^2}\right)\Delta' t^{-1}$$

and noting that

$$\frac{c_0}{a_1} = 2\nu_1 + t^{-1}\tilde{c}_5 \quad , \quad \frac{c_1}{a_1} = -2\nu_1 + t^{-1}\tilde{c}_6 \quad , \quad \tilde{c}_{5,6} \in \mathcal{F}$$

we separate the dominant linear terms and write

$$F_{1d} = 2\nu_1 t^{-2}\Delta - 2\nu_1 t^{-1}\Delta' \tag{3.31}$$

$$F_{1s} = t^{-3}\tilde{c}_5\Delta + t^{-2}\tilde{c}_6\Delta' \tag{3.32}$$

Finally, the terms which are at least quadratic are

$$\tilde{R}_2 = \frac{N}{a_1\Delta' + b_1} \tag{3.33}$$

where

$$N = \frac{q_3 {S_1}^2 {S_2}^2}{T^6 t^4} + \frac{q_2 {S_2}^2 S_1}{T^4 t^3} + \frac{q_5 S_1 {S_2}^3}{T^6 t^4} + \frac{d_1 S_2 {S_1}^2}{T^4 t^3} + \left(-\frac{c_0 a_1}{b_1 t^3 T^4} + \frac{d_2}{T^2 t^2}\right) S_2 S_1$$

$$+\frac{q_6 {S_2}^4}{T^6 t^4} + \frac{q_4 {S_2}^3}{T^4 t^3} + \frac{d_3 {S_1}^2}{T^2 t^2} + \left(-\frac{c_1 a_1}{b_1 t^3 T^4} + \frac{c_2 {a_1}^2}{t^5 T^6 {b_1}^2} + \frac{q_1}{T^2 t^2}\right) {S_2}^2 \tag{3.34}$$

3.6. **Estimate of $\tilde{R}_2(\Delta^{[1]}) - \tilde{R}_2(\Delta^{[2]})$.** The estimate is straightforward; below we provide details. Denote, for simplicity, $N(\Delta^{[j]}) = N^{[j]}$, $S_k(\Delta^{[j]}) = S_k^{[j]}$, $S^{[j]} = \max\{|S_1|^{[j]}, |S_2|^{[j]}\}$, $S = \max\{S^{[1]}, S^{[2]}\}$, $|\Delta| = \max\{|\Delta_1|, |\Delta_2|\}$.
We write

$$\left|\tilde{R}_2(\Delta^{[1]}) - \tilde{R}_2(\Delta^{[2]})\right| \le \frac{|N^{[1]} - N^{[2]}|}{|a_1\Delta_2^{[1]} + b_1|} + \frac{|N^{[2]}||a_1||\Delta_2^{[1]} - \Delta_2^{[2]}|}{|a_1\Delta_2^{[1]} + b_1||a_1\Delta_2^{[2]} + b_1|} \tag{3.35}$$

We have

$$|N^{[1]} - N^{[2]}| \le const|S^{[1]} - S^{[2]}|\left(t^{-2}T^{-2}S + t^{-3}T^{-4}S^2 + t^{-4}T^{-6}S^3\right)$$

(where the constant depends on $\rho_0$) and since $\|\Delta\| = \sup_t |tS|$ we get

$$|N^{[1]} - N^{[2]}| \le const\ t^{-4}T^{-2}\|\Delta^{[1]} - \Delta^{[2]}\| \tag{3.36}$$

Also

$$|N^{[2]}| \le const\ t^{-2}T^{-2}S^2 + t^{-3}T^{-4}S^3 + t^{-4}T^{-6}S^4 \tag{3.37}$$

The estimate (2.24) follows from (3.35), (3.36), (3.37).

3.7. **Estimate of $\tau$.** A direct calculation shows that (see (2.11) for the definition of $\tau$)

$$\tau = Q^2 F \tag{3.38}$$

where

$$F = \frac{a_2 - \Delta'}{(a_1 + b_1\Delta')^3} \frac{1}{\left(1 + \sqrt{1 - \frac{4T_0T_2}{T_1^2}}\right)^2}$$

and

$$Q = \frac{q_4 {S_2}^3}{T^3t^2} + \frac{q_3 {S_1}^2 {S_2}^2}{T^5t^3} + \frac{q_6 {S_2}^4}{T^5t^3} + \frac{d_3 {S_1}^2}{Tt} + \frac{c_1S_2}{Tt} + \frac{c_2}{Tt^2} + \frac{c_0S_1}{Tt}$$

$$+2\frac{d_1S_2{S_1}^2}{T^3t^2} + \frac{q_2{S_2}^2S_1}{T^3t^2} + \frac{d_2S_2S_1}{Tt} + \frac{q_1{S_2}^2}{Tt} \tag{3.39}$$

Note that on $\mathcal{B}_M$ we have

(3.40)
$|Q| \le KT^{-1}t^{-2}A(M,t)$ , where $A(M,t) = 1 + M + M^2t^{-1} + (M^3 + M^4)T^{-2}t^{-3}$

and

(3.41) $|F| \le KB(M,t)$ , where $B(M,t) = (1 + MT^{-2}t^{-2})(1 - MT^{-2}t^{-2})^{-3}$

and $K$ is a constant independent of $M$.

We use the notations of §3.6.

To estimate the difference $\tau^{[1]} - \tau^{[2]}$ of values of $\tau$ on two elements $\Delta^{[1]}, \Delta^{[2]}$ of $\mathcal{B}_M$ we write

$$|(Q^{[1]})^2F^{[1]} - (Q^{[2]})^2F^{[2]}| = |(Q^{[1]} - Q^{[2]})(Q^{[1]} + Q^{[2]})F^{[1]} + (Q^{[2]})^2(F^{[1]} - F^{[2]})|$$

$$\le 2QF\ |Q^{[1]} - Q^{[2]}| + Q^2|F^{[1]} - F^{[2]}|$$

Since

$$|Q| \le const\ \left(t^{-2}T^{-1} + t^{-1}T^{-1}(S + S^2) + t^{-2}T^{-3}S^3 + +t^{-3}T^{-5}S^4\right)$$

then on $\mathcal{B}_M$ we have $|Q| \le const\ t^{-2}T^{-1}$. Also

$$|F| \le const \tag{3.42}$$

Similarly

$$|Q^{[1]} - Q^{[2]}| \le const\ t^{-2}T^{-1}\|\Delta^{[1]} - \Delta^{[2]}\|$$

and

$$|F^{[1]} - F^{[2]}| \le const\left(|\Delta_2^{[1]} - \Delta_2^{[2]}| + t^{-1}|\Delta_1^{[1]} - \Delta_1^{[2]}|\right) \le t^{-3}T^{-2}\|\Delta^{[1]} - \Delta^{[2]}\|$$

The estimate (2.25) follows.

3.8. **Estimates under the assumptions of §2.5.** Note that in the assumptions of §2.5 we have $|S_{1,2}| \le \epsilon$.

3.8.1. *Estimate of $N$.* We split (3.34) in the form

$$N = N_{dom} + \delta_N$$

where

$$N_{dom} = \frac{d_2S_2S_1 + d_3{S_1}^2 + q_1{S_2}^2}{t^2T^2} = \frac{(d_2S_2 + d_3S_1)T^2\Delta + q_1S_2tT^2\Delta'}{t^2T^2}$$

and

$$\delta_N = \frac{c_2{a_1}^2{S_2}^2}{t^5T^6{b_1}^2} + \frac{S_2\left(-c_0S_1a_1 + q_4{S_2}^2b_1 + q_2S_2S_1b_1 + d_1{S_1}^2b_1 - c1a_1S_2\right)}{t^3b_1T^4}$$
$$+ \frac{{S_2}^2\left(q_5S_2S_1 + q_6{S_2}^2 + q_3{S_1}^2\right)}{T^6t^4}$$

Note that

$$|N_{dom}| \le const\ \epsilon\left(|t^{-2}\Delta| + |t^{-1}\Delta'|\right)$$

and

$$|\delta_N| \le const\ \epsilon^2|t^{-3}T^{-4}|$$

3.8.2. *Estimate of $Q$.* We split (3.39) in the form

$$Q = \frac{d_3{S_1}^2 + c_1S_2 + c_0S_1 + d_2S_2S_1 + q_1{S_2}^2}{tT} + \frac{c_2}{t^2T} + \delta \tag{3.43}$$

where

$$\delta = \frac{S_2\left(q_4{S_2}^2 + q_2S_2S_1 + d_1{S_1}^2\right)}{T^3t^2} + \frac{{S_2}^2\left(q_6{S_2}^2 + q_5S_2S_1 + q_3{S_1}^2\right)}{T^5t^3}$$

Hence from (3.43) we get

$$Q^2 = Q_{dom} + \delta_Q$$

where

$$Q_{dom} = \frac{1}{t^2T^2}\left[{d_3}^2{S_1}^4 + 2d_3\left(d_2S_2 + c_0\right){S_1}^3\right.$$
$$+ \left({d_2}^2{S_2}^2 + 2c_0d_2S_2 + {c_0}^2 + 2d_3q_1{S_2}^2 + 2d_3c_1S_2\right){S_1}^2$$
$$\left.+2S_2\left(d_2S_2 + c_0\right)\left(q_1S_2 + c_1\right)S_1 + {q_1}^2{S_2}^4 + {c_1}^2{S_2}^2 + 2c_1{S_2}^3q_1\right]$$

and

$$\delta_Q = \frac{{c_2}^2}{t^4T^2} + \frac{\delta^2}{t^4T^6} + 2\frac{\delta\left(q_1{S_2}^2 + c_0S_1 + d_2S_2S_1 + d_3{S_1}^2 + c_1S_2\right)}{t^3T^4}$$
$$+2\frac{c_2\left(q_1{S_2}^2 + c_0S_1 + d_2S_2S_1 + d_3{S_1}^2 + c_1S_2\right)}{T^2t^3} + 2\frac{c_2\delta}{t^4T^4}$$

Since we can rewrite

$$Q_{dom} = \frac{1}{t^2T^2}\left\{\left[{d_3}^2{S_1}^3 + 2d_3\left(d_2S_2 + c_0\right){S_1}^2\right.\right.$$
$$+ \left({d_2}^2{S_2}^2 + 2c_0d_2S_2 + {c_0}^2 + 2d_3q_1{S_2}^2 + 2d_3c_1S_2\right)S_1$$
$$\left.\left.+2S_2\left(d_2S_2 + c_0\right)\left(q_1S_2 + c_1\right)\right]T^2\Delta + \left({q_1}^2{S_2}^3 + {c_1}^2S_2 + 2c_1{S_2}^2q_1\right)tT^2\Delta'\right\}$$

we see that

$$|Q_{dom}| \leq const\ \epsilon \left(|t^{-2}\Delta| + |t^{-1}\Delta'|\right)$$

Also, clearly

$$|\delta_Q| \leq const\ \epsilon^2 |t^{-2}T^{-3}|$$

**Acknowledgments.** The authors are indebted to Professors P Deift and A Borodin for suggesting the problem and for many useful comments. The authors also acknowledge NSF support from grants 0100495 (OC) and 0074924 (RDC).